\documentclass[12pt]{article}
\usepackage[top=1in,bottom=1in,left=1in,right=1in]{geometry}
\usepackage{indentfirst}
\usepackage{amsfonts,amsmath,amsthm,amssymb}
\usepackage{longtable,enumerate,xcolor,seqsplit}
\usepackage{hyperref}

\theoremstyle{definition}
\newtheorem{definition}{Definition}
\newtheorem{remark}[definition]{Remark}
\newtheorem{question}{Question}
\newtheorem{problem}[question]{Problem}
\theoremstyle{plain}

\newtheorem{theorem}[definition]{Theorem}
\newtheorem{proposition}[definition]{Proposition}
\newtheorem{lemma}[definition]{Lemma}

\newtheorem{construction}[definition]{Construction}

\renewcommand{\leq}{\leqslant}
\renewcommand{\geq}{\geqslant}

\newcommand{\B}{\mathcal{B}}
\def\P{\mathcal{P}}

\newcommand{\PG}{{\rm PG}}
\newcommand{\AG}{{\rm AG}}
\def\proof{\par\noindent{\sc Proof~}}
\newcommand{\Aut}{{\rm Aut}}

\newcommand{\PSL}{{\rm PSL}}
\newcommand{\PGL}{{\rm PGL}}

\newcommand{\PGGL}{{\rm P\Gamma L}}

\newcommand{\AGL}{{\rm AGL}}

\newcommand{\CC}{\mathcal{C}}

\newcommand{\la}{\lambda}
\newcommand{\al}{\alpha}
\DeclareMathOperator{\Sym}{S} 
\DeclareMathOperator{\Alt}{A} 

\DeclareMathOperator{\D}{\mathcal{D}}

\providecommand{\keywords}[1]
{
  \small	
  \textbf{\textit{Keywords: }} #1
}

\begin{document}

\title{On flag-transitive imprimitive $2$-designs}

\author{Alice Devillers\thanks{This work was supported by an Australian Research Council Discovery Grant Project DP200100080.} and Cheryl E. Praeger\thanks{Address: Centre for the Mathematics of Symmetry and Computation, University of Western Australia, Perth, WA 6009, Australia.\quad  Emails: alice.devillers@uwa.edu.au, cheryl.praeger@uwa.edu.au}}
\maketitle

\begin{abstract}
In 1987, Huw Davies proved that, for a flag-transitive point-imprimitive $2$-$(v,k,\la)$ design, both the block-size $k$ and the number $v$ of points are bounded by functions of $\lambda$, but he did not make these bounds explicit. In this paper we derive explicit polynomial functions of $\lambda$ bounding $k$ and $v$. For $\lambda\leq 4$ we obtain a list of `numerically feasible' parameter sets $v, k, \lambda$ together with the number of parts and part-size of an invariant point-partition and the size of a nontrivial block-part intersection. Moreover from these parameter sets we determine all examples with fewer than $100$ points. There are exactly eleven such examples, and for one of these designs, a flag-regular, point-imprimitive $2-(36,8,4)$ design with automorphism group $\Sym_6$, there seems to be no construction previously available in the literature. 

\end{abstract}

\section{Introduction}

A \emph{$2$-$(v,k,\lambda)$ design} $\D$ is a pair $(\P,\B)$ with a set $\P$ of $v$ \emph{points} and a set $\B$ of \emph{blocks}
such that each block is a $k$-subset of $\P$ and each pair of distinct points is contained in $\lambda$ blocks.
We say $\D$ is \emph{nontrivial} if $2<k<v$, and \emph{symmetric} if $v=|\B|$.
All $2$-$(v,k,\lambda)$ designs in this paper are assumed to be nontrivial.
An automorphism of $\D$ is a permutation of the point set $\P$
which preserves the block set.
The set of all automorphisms of $\D$ under composition of permutations forms a
group, denoted by $\Aut(\D)$.
A subgroup $G$ of $\Aut(\D)$ leaves invariant a partition $\CC$ of $\P$ if each element
of $G$ permutes the parts of $\CC$ setwise. A partition $\CC$ is trivial 
if either $\CC$ consists of singleton sets, or $\CC=\{\P\}$; and $G$ is 
\emph{point-primitive} if the only $G$-invariant partitions of $\P$ are the 
trivial ones. Otherwise $G$ is said to be \emph{point-imprimitive}.
A \emph{flag} of $\D$ is a 
pair $(\alpha,B)$ where $\alpha\in\P$, $B\in\B$, and $B$ contains $\alpha$.
A subgroup $G$ of $\Aut(\D)$ is said to be \emph{flag-transitive} if $G$ acts transitively on the set of flags of $\D$.

A seminal result of Higman and McLaughlin \cite{HM} in 1961 showed that, in the case where $\lambda=1$, a flag-transitive subgroup of automorphisms is point-primitive. This break-through spurred others to discover whether this implication might hold more generally.
In particular Dembowski \cite[2.3.7(a)]{Dem} proved (in his 1968 book) that the same conclusion holds if $\lambda$ is coprime to the number $r$ of blocks containing a given point. However it does not hold in general.  
As pointed out by Davies \cite{Dav1}, Cameron and Kantor \cite[Theorem III]{CK79} showed  that the design whose points are the $2^{n+1}-1$ points of the projective space  $\PG(n,2)$, $n$ odd, and whose blocks are the hyperplane-complements, with natural incidence, admits $ \PGGL((n+1)/2,4)$ as a flag-transitive group of automorphisms that is point-imprimitive. For these designs $\la = 2^{n-1}$ grows exponentially with $n$.  
 
On the other hand, in  \cite{Dav1}, Davies also
established that, for fixed $\lambda$, there are only finitely many flag-transitive, point-imprimitive 2-designs, by showing that the block-size $k$ and the number $v$ of points are both bounded in terms of $\la$. However he did not give explicit upper bounds.  Some years later Cameron and the second author \cite[Proposition 4.1]{CP93} showed that also $v\leq (k-2)^2$, that is, $v$ is bounded above in terms of $k$, for flag-transitive, point-imprimitive designs, and that 
the smallest possible block size is $k= 6$. Recently Zhan and Zhou \cite{ZZ18} found that there are exactly 14 examples with $k=6$,  all with $v= (k-2)^2=16$.

Davies' examples above from projective geometry are all \emph{symmetric} designs, and indeed much progress has been made studying flag-transitive symmetric $2$-designs.
In \cite{Reg05},  O'Reilly-Regueiro showed that a  flag-transitive, point-imprimitive, symmetric design must have $k\leq\lambda(\lambda+1)$, and further work (see  \cite{LPR09, P07, PZ06}) refined this bound and classified all examples with $\lambda$ up to $4$. 

In this paper we find explicit bounds for $k$, and hence for $v$, in terms of $\la$, without assuming that the design is symmetric.

\begin{theorem}\label{main}
 Let $\cal D =(\P, \cal B)$ be a $2$-$(v,k,\la)$ non-trivial design admitting 
a flag-transitive point-imprimitive group  of automorphisms. Then $k\leq 2\la^2(\la-1)$ and $v\leq \left(2\la^2(\la-1)-2\right)^2$.
\end{theorem}

Recall the result of Higman and McLaughlin \cite{HM} that for $\lambda=1$, all flag-transitive $2$-designs are point-primitive. A recent result of the authors and colleagues in \cite[Theorem 1.1]{DLPX} shows that there are up to isomorphism just two designs which prevent this conclusion holding also for $\lambda=2$: namely there are exactly two  flag-transitive, point-imprimitive $2$-$(v,k,2)$ designs, and both of them are $2-(16,6,2)$ designs.
For the cases $\lambda=3,4$, we  list in Proposition~\ref{la=3} all `numerically feasible' parameter sets for flag-transitive, point-imprimitive $2$-designs, that is to say, parameter sets which satisfy all the conditions imposed by our preliminary results in Section~\ref{prelim}.
We specify not only the parameters $v, k, \lambda$, but also the number $d$ of parts and the part-size $c$ of a nontrivial invariant point-partition, and the (constant) size $\ell$ of a non-empty intersection between a part of this partition and a block of the design (see Lemma~\ref{le}). Although we have not managed to complete the classification of all examples with $\lambda\leq 4$, which is given in \cite{LPR09, P07, PZ06} in the symmetric case, we have been able to classify all examples with less than $100$ points, and in so-doing, we constructed a design on $36$ points for which the full automorphism group is flag-regular (and so no proper subgroup is flag-transitive). We thought this design, which has the parameters in line 5 of Table~$\ref{tableintro}$, was new (we asked a few design experts and none had seen it before), but after finishing our analysis we discovered that the design was identified by Zhang and Zhou in \cite[Theorem 1.3]{ZZ}. However  no construction of the design is given in \cite{ZZ}, see Remark~\ref{rem:shenglin}. 
In Section~\ref{sec:36}, we give several constructions and discuss this interesting design further.

\begin{theorem}\label{class}
 There are exactly eleven $2$-$(v,k,\la)$ non-trivial designs admitting 
a flag-transitive point-imprimitive group $G$ of automorphisms, with $\la\leq 4$ and $v<100$, with two of them admitting two partitions of different sizes. If $G$ preserves a partition into $d$ parts of size $c$, then $(\lambda,v, k,r,c,d)$ are as in one of the lines of Table $\ref{tableintro}$, the penultimate column gives the number of designs with these parameters (up to isomorphism), and the last column gives a reference when possible. 

\end{theorem}
\begin{table}[h]
\centering
\caption{\label{tableintro}The eleven designs for $\la\leq 4$ and $v<100$}
 \begin{tabular}{|cccccc|c|c|}
\hline
$\la$&$v$&$k$&$r$&$c$&$d$&Number &Reference\\
\hline
$2$&$16$&$6$&$6$&$4$&$4$&$2$&{\rm \cite{AS, Burau, Huss}}\\
$3$&$45$&$12$&$12$&$9$&$5$&$1$&{\rm \cite{P07}}\\
$4$&$15$&$8$&$8$&$3$&$5$&$1$&{\rm \cite{CK79, Dav1}}\\
$4$&$16$&$6$&$12$&$4$&$4$&$2$&{\rm \cite{ZZ18}}\\
$4$&$36$&$8$&$20$&$6$&$6$&$1$&Construction~{\rm\ref{con1}}, \cite{ZZ}\\
$4$&$96$&$20$&$20$&$6$&$16$&$2$&{\rm \cite{LPR09}}\\
$4$&$96$&$20$&$20$&$16$&$6$&$4$&{\rm \cite{LPR09}}\\
\hline
\end{tabular}
   
\end{table}

\begin{remark}
(a)  We note that the two designs with $(\lambda,v, k,r,c,d)=(4,96,20,20,6,16)$  in Table~\ref{tableintro} 
are among the four flag-transitive designs for $(4,96,20,20,16,6)$ (see the classification in \cite{LPR09}). Thus there are exactly eleven  designs satisfying the conditions of Theorem~\ref{class}, two of which admit two nontrivial partitions with different parameters $(c,d)$. See Remark~\ref{rem:t2} for more details.

(b)
The smallest value of $\lambda$ for which flag-transitive, point-imprimitive designs may exist is $\lambda=2$ and, as we mentioned above, in this case it follows from \cite[Theorem 1.1]{DLPX} that $\cal D$ is one of two known $2-(16,6,2)$ designs. Thus the upper bounds on $(k, v)$ in Theorem~\ref{main} when $\lambda=2$, namely $(8,36)$, are far from tight. Also, for $\lambda=4$, it follows from  Proposition~\ref{la=3}  that the bounds on both $k$ and $v$ in Theorem~\ref{main} are definitely not tight. If $\lambda=3$, 
then the value of $k$ could possibly meet the bound $k=36$ of Theorem~\ref{main}, with the remaining parameters as in one of three lines of the table in Proposition~\ref{la=3}.  
Thus we ask in general:
\end{remark}

\begin{question}\label{q1}
Can the functions of $\lambda$ bounding $k$ and $v$ in Theorem~\ref{main} be improved?
\end{question}

We think the answer to Question~\ref{q1} is `yes' (with the possible exception of $\lambda=3$) and would like to see improved polynomial bounds.
If $\lambda=3$, then an answer to the next question would settle the tightness of the bounds in Theorem~\ref{main} for that case.

\begin{question}\label{q2}
Does there exist a flag-transitive, point-imprimitive design with parameter set 
\[
(\lambda,v, k,r,c,d) = (3,561, 36, 48, 17, 33),\  (3,561, 36, 48, 33, 17), \text{ or } 
(3, 1156, 36, 99, 34, 34)?
\]
\end{question}

When $\lambda=3$, there are seven lines of the table in Proposition~\ref{la=3} which have not been treated in Theorem~\ref{class}, that is to say, four lines in addition to the parameter sets in Question~\ref{q2}.  Also, for $\lambda=4$, there are eleven lines of the table in Proposition~\ref{la=3} which have not been treated in Theorem~\ref{class}.

\begin{problem}\label{prob}

Classify all the flag-transitive, point-imprimitive $2-(v,k,\lambda)$ designs with parameter sets $(\lambda,v, k,r,c,d)$ as in one of the 18 lines of the tables in Proposition~\ref{la=3} with $v\geq 100$.

\end{problem}

A complete answer to Problem~\ref{prob} would finish the classification of the flag-transitive, point-imprimitive $2-(v,k,\lambda)$ designs with $\lambda\leq 4$. 
A partial answer to Problem~\ref{prob} is given in \cite[Theorem 1.3]{ZZ} under the additional assumption that the flag-transitive, point-imprimitive group is also point-quasiprimitive (that is, all nontrivial normal subgroups are point-transitive). Thus when attacking Problem~\ref{prob}, one may assume that the group is not point-quasiprimitive.

In Section \ref{prelim}  we list some well-known facts about $2$-designs and  prove some numerical conditions for flag-transitive point-imprimitive $2$-designs.
In Section \ref{sec:main} we prove Theorem \ref{main}. In Section \ref{sec:small} we determine all numerically feasible parameters sets for $\la=3, 4$. In Section \ref{sec:36} we give several constructions for a $2$-design on $36$ points, and we show that up to isomorphism this design is the unique flag-transitive, point-imprimitive $2-(36,8,4)$ design (Proposition~\ref{lem:36unique}). 
Finally in Section \ref{sec:class}, we classify all flag-transitive point-imprimitive $2$-designs with  $\la\leq 4$ and $v<100$, providing lots of information on their automorphism groups and how to construct them with Magma \cite{magma}.

\section{Preliminary results on designs}\label{prelim}
We first collect some useful results on flag-transitive designs.

\begin{lemma} \label{condition 1}
Let $\D =(\P, \cal B)$ be a $2$-$(v,k,\lambda)$ design and let $b=|\cal B|$.
Then the number of blocks of $\D$ containing each point of $\D$ is a constant $r$ satisfying the following:
\begin{enumerate}
\item[\rm(i)] $r(k-1)=\lambda(v-1)$;
\item[\rm(ii)] $bk=vr$;
\item[\rm(iii)] $b\geq v$ and $r\geq k$;
\item[\rm(iv)] $r^2>\lambda v$.
\end{enumerate}
In particular, if $\D$ is not symmetric then $b>v$ and $r>k$.
\end{lemma}

\proof
Parts~(i) and~(ii) follow immediately by simple counting.
Part~(iii) is  Fisher's Inequality \cite[p.99]{Ryser}.
By~(i) and~(iii) we have
\[
r(r-1)\geq r(k-1)=\lambda(v-1)
\]
and so $r^2\geq\lambda v+r-\lambda$.
Since $\D$ is nontrivial, we deduce from (i) that $r>\lambda$.
Hence $r^2>\lambda v$, as stated in part~(iv).

\qed

We now prove the following important technical proposition.

\begin{lemma}
\label{le} Let $\cal D =(\P, \cal B)$ be a nontrivial $2$-$(v,k,\la)$ design admitting 
a flag-transitive point-imprimitive group $G$ of automorphisms, which leaves invariant a nontrivial point-partition $\mathcal{C}$ into $d$ parts of size $c$. Then the non-empty intersections $B\cap \Delta$, for $B\in\cal B$ and $\Delta\in\mathcal{C}$, have a constant size $\ell$, say. Moreover, the integer $x=k-1-d(\ell-1)$ is  positive, and the following equalities, inequalities and divisibility conditions hold:
\begin{enumerate}[(i)]
\item $\la\geq 2$; 
\item $\ell\mid k\quad \mbox{and}\ 1<\ell < k$;
\item $\la(c-1)=r(\ell-1)$;
\item $k=xc+\ell$;
\item $rx=\la (d-1)$;
\item $k\mid \frac{\la c(c-1)(k-(x+1))}{(\ell-1)^2}$, in particular, if $\ell=2$ then $k\mid \la c(c-1)(x+1)$;
\item $k\mid \la\ell(x+1)(x+\ell)$;
\item $x(\ell-1)\leq \la-1$;
\item $c\geq \frac{\la+\ell(\ell-1)}{\la-x(\ell-1)}$;
\item $k\geq \frac{\la(x+\ell)}{\la-x(\ell-1)}$;
\end{enumerate} 
\end{lemma}

\proof
By the celebrated result of Higman and McLaughlin \cite{HM} mentioned above, if a $2-(v,k,1)$ design (linear space) is flag-transitive, then it is point-primitive. Thus $\la\geq 2$, proving (i).
Let 
$\mathcal{C}=\{\Delta_1,\Delta_2,\dots,\Delta_d\}$,  with $1<d<v$ and $|\Delta_i|=c>1$ for each $i$,
so that
\begin{equation}\label{Eq1}
v=cd.
\end{equation}

Let $B, B'\in \mathcal{B}$ and $\Delta, \Delta'\in \mathcal{C}$ such that $B\cap\Delta$ and $B'\cap\Delta'$ are non-empty, and choose $\alpha\in B\cap\Delta$ and $\alpha'\in B'\cap\Delta'$. Since $G$ is flag-transitive, there exists $g\in G$ such that $(B,\alpha)^g=(B',\alpha')$. As $\alpha^g=\alpha'$ we have $\Delta^g=\Delta'$, and hence  $(B\cap\Delta)^g= B'\cap\Delta'$. Thus $\ell=|B\cap\Delta|$ is independent of $B$ and $\Delta$, and so $\ell\mid k$. Since $G$ is block-transitive and $\cal D$ is a $2$-design, it follows that each block contains a pair of points in the same part of $\mathcal{C}$, and a pair of points from different parts of $\mathcal{C}$. Thus $ 1<\ell < k$, and this proves (ii).

Fix a point $\alpha$, a block $B$  containing $\alpha$, and let $\Delta$ be the part of $\mathcal{C}$ containing $\alpha$.
Counting the point-block pairs $(\alpha',B')$ with $\alpha'\in\Delta\setminus\{\alpha\}$ and $B'$ containing $\alpha$ and $\alpha'$, we obtain
$
\la(c-1)=r(\ell-1),
$ 
proving (iii).
Multiplying both sides of this equation by $k-1$, and using Lemma~\ref{condition 1}(i) and equation \eqref{Eq1}, we find that $\lambda(k-1)(c-1)=r(k-1)(\ell-1) = \lambda (v-1)(\ell-1)$, and hence that
$$
(cd-1)(\ell-1)=(k-1)(c-1).
$$
Thus $cd(\ell-1)-(\ell-1)=c(k-1)-(k-1)$, from which we deduce that $k-\ell=c\left(k-1-d(\ell-1)\right)$.
Since $x=k-1-d(\ell-1)$ and since $\ell<k$, this implies that $x$ is a positive integer. Also it follows from this equation that 
$
k=xc+\ell,$ proving (iv).

Using Lemma~\ref{condition 1}(i), part (iii) and \eqref{Eq1}, we get that 
$$
rx=r(k-1)-dr(\ell-1)=\la(v-1)-d\la(c-1)=\la (d-1),
$$ proving (v).
%

%
By  part (v) we have $d= 1+(rx/\lambda)$, and by (iii), $r = \lambda(c-1)/(\ell-1)$, so that $d=1+ x(c-1)/(\ell-1)$.Then part (iv) and \eqref{Eq1} imply that  
\begin{align*}
vr&=cdr=c\left(\frac{c-1}{\ell-1}\cdot x+1\right)\left( \frac{\la(c-1)}{\ell-1}\right)\\
&=\frac{\la c(c-1)(cx-x+\ell-1)}{(\ell-1)^2}=\frac{\la c(c-1)(k-(x+1))}{(\ell-1)^2}.
\end{align*}

By Lemma~\ref{condition 1}(ii), $bk=vr$, so  
$$
k\mid \frac{\la c(c-1)(k-(x+1))}{(\ell-1)^2}.
$$ In particular, for $\ell=2$, $k\mid\la c(c-1)(k-(x+1))$ and thus also $k\mid\la c(c-1)(x+1)$, proving (vi).
It follows that 
$k$ divides $\la c(c-1)(x+1)$, and hence $k$ also divides $\la (xc)(xc-x)(x+1)$, which by part (iv) is equal to  $\la (k-\ell)(k-\ell-x)(x+1)$. Thus
$$k\mid \la \ell(\ell+x)(x+1),$$ proving (vii).

On the other hand, since $r\geq k$ (by Lemma~\ref{condition 1}(iii)), and using part (iv) and part (iii), we have
\[
\la k>\la(k-\ell-x)=\la x(c-1)=rx(\ell-1)\geq kx(\ell-1),
\]
and so $x(\ell-1)<\la$. Since all these parameters are integers, (viii) follows.

Now using part (iii), the inequality $r\geq k$, and 
part (iv), we find 
$$
\la(c-1)=r(\ell-1) \geq k(\ell-1) = (xc+\ell)(\ell-1).
$$ 
Rearranging this inequality gives $c(\lambda - x(\ell-1))\geq \lambda +\ell(\ell-1)$, and since  $\la-x(\ell-1)>0$ by part (viii), we have 
$$
c\geq\frac{\la+\ell(\ell-1)}{\la-x(\ell-1)},
$$ 
proving (ix).  Finally, using (iv) and (ix), 
$$
k=xc+\ell\geq x.\frac{\la+\ell(\ell-1)}{\la-x(\ell-1)}+\ell=\frac{\la(x+\ell)}{\la-x(\ell-1)},
$$ 
proving (x).
\qed

\medskip
We will need the following technical lemma.

\begin{lemma}\label{le2} 
Let $z$ be a real number greater than $1$.
The function $g(x,y)=(x+1)(y+1)(x+y+1)$ from $\mathbb{R}^2$ to $\mathbb{R}$, restricted to the hyperbola $xy=z$ with $x,y\geq 1$ decreases as $x$ increases between $1$ and $\sqrt{z}$, increases as  $x$ increases between $\sqrt{z}$ and $z$, and has a maximum of $2(z+1)(z+2)$ at $(x,y)=(1,z)$ and $(z,1)$.
\end{lemma}

\proof
On the hyperbola $xy=z$, the function $g$ becomes 
\begin{align*}
g(x,z/x)&=(x+1) (z/x+1)(x+z/x+1)\\
&=\frac{(x+1)(x+z)(x^2+x+z)}{x^2} \\
&=\frac{x^4+(z+2)x^3+(3z+1)x^2+z(z+2)x+z^2}{x^2} \\
&=x^2+(z+2)x+(3z+1)+z(z+2)x^{-1}+z^2 x^{-2}.
\end{align*}

We can now compute the derivative 
\begin{align*}
g'(x,z/x)&=2x+(z+2)-z(z+2)x^{-2}-2z^2 x^{-3} \\
&=\frac{2x^4+(z+2)x^3-z(z+2)x-2z^2 }{x^3}\\
&=\frac{2(x^2-z)(x^2+z)+(z+2)x(x^2-z) }{x^3}\\
&=\frac{(x^2-z)(2(x^2+z)+(z+2)x) }{x^3}.
\end{align*}
Since $x\geq 1$ and $z>1$,
 the denominator and second factor of the numerator are obviously positive, while the first factor of the numerator is negative when $x<\sqrt{z}$ and positive when $x>\sqrt{z}$. Therefore the maximum of $g(x,y)$  on the hyperbola is  $g(1,z)=g(z,1)=2(z+1)(z+2)$.
\qed

\medskip

\section{Proof of Theorem~\ref{main}}\label{sec:main}
The preparatory results from Section \ref{prelim} allow us to obtain our first bound.

\begin{proposition}\label{firstbound}
 Let $\cal D =(\P, \cal B)$ be a $2$-$(v,k,\la)$ non-trivial design admitting 
a flag-transitive point-imprimitive group  of automorphisms. Then $k\leq 2\la^2(\la+1)$.
\end{proposition}
\proof Let 
all parameters be as in Lemma \ref{le}.
Then 
$
k\leq \la \ell(\ell+x)(x+1),
$ by Lemma \ref{le}(vii).
At this point, it is convenient to change variables: let $y=\ell-1$ and $\mu=\la-1$, so that $1\leq y\leq \mu$, $\mu\geq 1$ and $xy\leq \mu$ by Lemma \ref{le}(viii).
Let $$g(x,y)=(x+1)(y+1)(x+y+1),$$
so that $k\leq (\mu+1) g(x,y)$. We wish to find the maximum of the function $g(x,y)$ on the domain $x,y\geq 1$, $xy\leq \mu$.

Since $g(x,y)$ increases with $y$, for a fixed $x$, the maximum of this function must be on the hyperbola $xy=\mu$.
By Lemma \ref{le2}, the maximum of $g(x,y)$ on that hyperbola is $2(\mu+1)(\mu+2)$ obtained at $(x,y)=(1,\mu)$ and $(\mu,1)$.
 Therefore $k\leq 2(\mu+1)^2(\mu+2)= 2\la^2(\la+1)$. 
\qed

\medskip
For $\la=2$ the bound in Proposition~\ref{firstbound} gives $k\leq 24$. Together with Liang and 
Xia, the authors showed in \cite{DLPX} that there are only two imprimitive flag-transitive 2-designs, both of which are $2-(16,6,2)$ designs. Thus this bound is definitely not tight for all $\lambda$.

For $\la=3$ the bound in Proposition~\ref{firstbound} gives  $k\leq 72$. That is also not the best possible, as looking at $\la=3$ in detail  (splitting up into cases for possible $(\ell,x)$) we can show $k\leq 36$, see Proposition \ref{la=3} (note our list below matches with the cases listed in \cite{Dav1}).

For $\la=4$ the bound in Proposition~\ref{firstbound} gives $k\leq 160$ (which is better than what is stated in \cite{Dav1}), but we improve this in Proposition \ref{la=3} to $k\leq 80$. 

Now we prove the main theorem. Note that we use here Proposition \ref{la=3} which is in the next section. That proposition only relies on results from Section \ref{prelim} so our argument is not circular.

\medskip

\proof[Theorem \ref{main}] Let $G$ be the flag-transitive automorphism group, and let $\mathcal{C}$ be the $G$-invariant non-trivial partition.  Let 
all parameters be as in Lemma \ref{le}. 
If $\la=2$, we showed in \cite{DLPX} (by group-theoretic arguments) that $k=6<2\la^2(\la-1)$.
The statement is clearly true for $3\leq \la\leq 4$ by Proposition \ref{la=3} below, so assume $\la\geq 5$.

We first claim that $k$ cannot be equal to the bound found in Proposition \ref{firstbound}. Assume to the contrary that $k=2\la^2(\la+1)$. Looking at the proof of Proposition \ref{firstbound}, this implies that  $k=\la \ell(\ell+x)(x+1)$, $x(\ell-1)=\la-1$, and $\ell=2$ or $\la$. 
If $\ell=2$, then $x=\la-1$, and $x$ divides $k-\ell=2(\la^3+\la^2-1)$ by Lemma \ref{le}(iv). It follows that $\la-1$ divides $2$, so $\la=2$ or $3$, a contradiction since  $\la\geq 5$.
If $\ell=\la$, then $x=1$, $c=k-\ell=2\la^3+2\la^2-\la$ by Lemma \ref{le}(iv), and  $\la-1$ divides $\la(c-1)=\la(2\la^3+2\la^2-\la-1)$ by Lemma \ref{le}(iii). It also follows that $\la-1$ divides $2$, so $\la=2$ or $3$, contradicting  $\la\geq 5$.
Thus the claim is proved.

For convenience, we now use  the notation $y=\ell-1$ and $\mu=\la-1\geq 4$, as in the proof of Proposition \ref{firstbound}. Recall that $k$ divides $\la \ell(x+1)(x+\ell)=(\mu+1)g(x,y)$ by Lemma~\ref{le}(vii), and we have just shown that $k<2\lambda^2(\lambda+1)$.

We claim that $k\leq \max\{\lambda^2(\lambda+1), (\mu+1) X\}$, where $X$ is the  second largest value of $g(x,y)$ on the domain  $x,y\geq 1$, $xy\leq \mu$ with $x,y$ integers. We see this as follows. In Proposition \ref{firstbound}, it was shown that the maximum value of $g(x,y)$ on this domain is $2(\mu+1)(\mu+2)=2\lambda(\lambda+1)$.  If $g(x,y)$ takes this value, that is, if $(\mu+1)g(x,y)=2\la^2(\la+1)$,  then $k$ must be a proper divisor and hence $k\leq \lambda^2(\lambda+1)$. On the other hand, if 
$(\mu+1)g(x,y)<2\la^2(\la+1)$, then $k\leq (\mu+1)g(x,y)\leq (\mu+1) X$. This proves the claim.

We now determine $X$, the second largest value of $g(x,y)$. Note that, if $xy<\mu$, then by Lemma \ref{le2}, $g(x,y)\leq g(1,xy)\leq g(1,\mu-1)=2\mu(\mu+1)$.
Hence $X$ is either $2\mu(\mu+1)$, or  $g(x,\mu/x)$ for some integer $x\neq 1$ properly dividing $\mu$. If $\mu$ is prime then such an $x$ does not exist, and   $X=2\mu(\mu+1)$.

Assume $\mu$ is not a prime, and let $p$ be the smallest prime factor of $\mu$, so that $\mu=pq$ where $q\geq p$. Note that $p\leq \sqrt{\mu}$. It follows from Lemma \ref{le2} that the largest value for $g(x,\mu/x)$, for some integer $x\neq 1$ dividing $\mu$, is $g(p,q)=(p+1)(q+1)(p+q+1)$. 
We claim that $2\mu(\mu+1)>(p+1)(q+1)(p+q+1)$, with the unique exception of $\mu=4$.

Using that $\mu=pq$, the above inequality can be rewritten as  
$$
q^2(2p^2-p-1)-q(p^2+p+2)-(p+1)^2>0.
$$
Suppose first that $p$ is odd. Then $2p^2-p-1\geq p^2+p+2,$ and it is elementary (taking the derivative with respect to $q$) to see that the left-hand side grows with $q$, and so 
$$
q^2(2p^2-p-1)-q(p^2+p+2)-(p+1)^2\geq p^2(2p^2-p-1)-p(p^2+p+2)-(p+1)^2=2p^4-2p^3-3p^2-4p-1,
$$ 
which is positive, so the claimed inequality holds. Suppose now that $p=2$. If $p=2$ and $q=3$ then 
$$
q^2(2p^2-p-1)-q(p^2+p+2)-(p+1)^2=12,
$$ 
and so $q^2(2p^2-p-1)-q(p^2+p+2)-(p+1)^2>0$ for $p=2$ and $q\geq 3$, but 
$$
q^2(2p^2-p-1)-q(p^2+p+2)-(p+1)^2<0
$$ 
for $p=q=2$, that is for $\mu=\la-1=4.$ This proves the claim.

To summarise,
either $X=2\mu(\mu+1)$, or $\mu=4$ in which case $X=g(2,2)=45$.
Recall from above that  $k\leq \max\{\la^2(\la+1),(\mu+1)X\}$.

If $\la>5$, then $(\mu+1)X=2\mu(\mu+1)^2=2\la^2(\la-1),$ which is larger than  $\la^2(\la+1)$, so $k\leq 2\la^2(\la-1)$.
Assume now that $\la=5$. Then $(\mu+1)X=225,$ which is larger than  $\la^2(\la+1)=150$, so $k\leq 225$. We claim that $k\leq 2\la^2(\la-1)=200$. Assume for a contradiction that  $200<k\leq 225$. If $x$ or $y$ is equal to $1$ then $k$ divides $60,120,200$ or $300$ by Lemma \ref{le}(vii), but we have seen above that $k\neq 300$, so in all those cases we have $k\leq 200$. Thus we must have that $x=y=2$ in which case  $k$ divides $225$. This implies that  $k=225$. Using Lemma \ref{le} parts (iv), (iii), (v) to get $c$, $r$, $d$ respectively, we deduce that $(c,d,k,r,\ell)=(111,111,225,275,3)$. Let $\Delta\in\mathcal{C}$, $D=G^\mathcal{C}$, and $L=(G_{\Delta})^{\Delta}$. By \cite[Theorem 5.5]{PS} we may assume that $G\leq L\wr D\leq \Sym_{111}\wr \Sym_{111}$, acting imprimitively. Also, by the argument above, there cannot be any other values of $c,d$ yielding a flag-transitive example, and hence both  $L$ and $D$ are primitive of degree $111$. The only such groups are $\Alt_{111}$ and $\Sym_{111}$. 
Let $\al, \beta$ be distinct points of $\Delta$, and let $B_1, \dots, B_5\in\mathcal{B}$ be the $\lambda=5$ blocks containing $\{\al,\, \beta\}$. Then $G_{\{\al, \beta\}}$, which is a subgroup of $G_\Delta$, fixes $X:=\cup_{i=1}^5(B_i\cap\Delta)$ setwise, and since each $B_i\cap\Delta$ has size $3$, the set $X$ has size $s$ where $3\leq s\leq 7$. On the other hand, $(G_{\{\al,\beta\}})^\Delta$ is the setwise stabiliser in $L$ of the pair $\{\al,\beta\}$. Since $L$ is $\Alt_{111}$ or $\Sym(111)$, $G_{\{\al,\beta\}}$ has orbits in $\Delta$ of lengths $2, 109$, which is a contradiction.  
Therefore, for $\lambda=5$, we also have that $k\leq 2\la^2(\la-1)$.
This proves the upper bound for $k$ in all cases.  By \cite[Proposition 4.1]{CP93}, $v\leq (k-2)^2$, and the claimed upper bound on $v$ follows. 
This finishes the proof.
\qed

\section{Numerically feasible parameter tuples for small \texorpdfstring{$\lambda$}{lambda}}\label{sec:small}

Recall that by the theorem of Higman and McLaughlin~\cite{HM}, $\la\neq 1$. 
By \cite[Theorem 1.1]{DLPX}, when $\la=2$,  the only values for $v, k, r$ are $v=16$, $k=r=6$. Lemma \ref{le}(viii) yields $\ell=2$, and by  Lemma \ref{le}(iv), $c=k-2=4$ and hence $d=v/c=4$ also. For this reason, we only need to consider $\la\geq 3$.

For specific small $\la$, we can list all the pairs $(\ell,x)$ satisfying Lemma \ref{le}(viii), and do a more refined investigation leading to all the possible tuples $(\la,v,k,r,c,d,\ell)$ that are \emph{numerically feasible}, in the sense that they satisfy all of the restrictions from Lemmas~\ref{condition 1} and \ref{le}.

 \begin{proposition}\label{la=3}
Suppose that $\la\in\{3,4\}$.
Then the numerically feasible parameters $(\la,v,k,r,c,d,\ell)$ for  a $2$-$(v,k,\la)$ non-trivial design admitting 
a flag-transitive point-imprimitive group of automorphisms  are as in one of the rows of the following tables, where $c,d,\ell$ are as in Lemma \ref{le} and $r$ is the number of blocks through a point. 

\begin{minipage}{.5\linewidth}
\centering
\begin{tabular}{|ccccccc|}
\hline
$\la$&$v$&$k$&$r$&$c$&$d$&$\ell$\\
\cline{1-7}
$3$&$16$&$6$&$9$&$4$&$4$&$2$\\
$3$&$45$&$12$&$12$&$5$&$9$&$2$\\
$3$&$45$&$12$&$12$&$9$&$5$&$3$\\
$3$&$100$&$12$&$27$&$10$&$10$&$2$\\
$3$&$120$&$18$&$21$&$8$&$15$&$2$\\
$3$&$120$&$18$&$21$&$15$&$8$&$3$\\
$3$&$256$&$18$&$45$&$16$&$16$&$2$\\
$3$&$561$&$36$&$48$&$17$&$33$&$2$\\
$3$&$561$&$36$&$48$&$33$&$17$&$3$\\
$3$&$1156$&$36$&$99$&$34$&$34$&$2$\\
\hline
\end{tabular} 
\end{minipage}
\begin{minipage}{.5\linewidth}
\centering
\begin{tabular}{|ccccccc|}
\hline
$\la$&$v$&$k$&$r$&$c$&$d$&$\ell$\\
\hline
$4$&$15$&$8$&$8$&$3$&$5$&$2$\\
$4$&$16$&$6$&$12$&$4$&$4$&$2$\\
$4$&$36$&$8$&$20$&$6$&$6$&$2$\\
$4$&$45$&$12$&$16$&$5$&$9$&$2$\\
$4$&$45$&$12$&$16$&$9$&$5$&$3$\\
$4$&$96$&$20$&$20$&$6$&$16$&$2$\\
$4$&$96$&$20$&$20$&$16$&$6$&$4$\\
$4$&$100$&$12$&$36$&$10$&$10$&$2$\\
$4$&$196$&$16$&$42$&$14$&$14$&$2$\\
$4$&$231$&$24$&$40$&$11$&$21$&$2$\\
$4$&$231$&$24$&$40$&$21$&$11$&$3$\\
$4$&$280$&$32$&$36$&$10$&$28$&$2$\\
$4$&$280$&$32$&$36$&$28$&$10$&$4$\\
$4$&$435$&$32$&$42$&$15$&$29$&$2$\\
$4$&$484$&$24$&$84$&$22$&$22$&$2$\\
$4$&$1976$&$80$&$100$&$26$&$76$&$2$\\
$4$&$1976$&$80$&$100$&$76$&$26$&$4$\\
$4$&$2116$&$48$&$180$&$46$&$46$&$2$\\
\hline
\end{tabular}
\end{minipage}

\end{proposition}
Note that in each case the number of blocks can be determined using the formula given  by Lemma~\ref{condition 1}(ii).
\medskip
\proof It is easy to check that all parameter sets in the tables satisfy all conditions from Lemmas~\ref{condition 1} and \ref{le}, with $x=k-1-d(\ell-1)$.
We now show that for each $\la=3,4$, the parameters must be as in one of the tables.



\medskip\noindent
\fbox{The case $\la=3$}\quad
Lemma \ref{le}(viii) yields  three possibilities for $(\ell,x)$, namely  $(2,1)$,  $(2,2)$, $(3,1)$.
We now split the analysis into these 3 cases.
\begin{enumerate}[(i)]
\item $(\ell,x)=(2,1)$. By Lemma \ref{le}(vii) $k\mid 36$ and by Lemma \ref{le}(x) $k\geq {9/2}$. Moreover $k$ is even by Lemma \ref{le}(ii).
Thus $k\in\{6,12,18,36\}$. By Lemma \ref{le}(iv) $c=k-2$, and Lemma \ref{le}(vi) yields $k\mid {6 c(c-1)}$ which is satisfied in each case. Combining Lemma \ref{le}(iii) and (v) yields $d=c$. By Lemma \ref{le}(iii) $r=3(c-1)$.
Thus the possibilities for $(c,d,k,r,\ell)$ are $(4,4,6,9,2)$,  $(10,10,12,27,2)$, $(16,16,18,45,2)$, $(34,34,36,99,2)$.

\item $(\ell,x)=(2,2)$. 
By Lemma \ref{le}(vii) $k\mid 72$ and by Lemma \ref{le}(x) $k\geq 12$. 
Thus $k\in\{12,18,24,36,72\}$. By Lemma \ref{le}(iv) $c=k/2-1$, and Lemma \ref{le}(vi) yields $k\mid {9 c(c-1)}$ which is not satisfied for $k=24$ and $k=72$. Combining Lemma \ref{le}(iii) and (v) yields $d=2c-1$. By Lemma \ref{le}(iii) $r=3(c-1)$.
Thus the possibilities for $(c,d,k,r,\ell)$ are $(5,9,12,12,2)$, $(8,15,18,21,2)$, $(17,33,36,48,2)$.

\item $(\ell,x)=(3,1)$. 
By Lemma \ref{le}(vii) $k\mid 72$ and by Lemma \ref{le}(x) $k\geq 12$. 
Thus $k\in\{12,18,24,36,72\}$. By Lemma \ref{le}(iv) $c=k-3$, and Lemma \ref{le}(vi) yields $k\mid \frac{3 c(c-1)(c+1)}{4}$ which is not satisfied for $k=24$ and $k=72$. Combining Lemma \ref{le}(iii) and (v) yields $d=(c+1)/2$.  By Lemma \ref{le}(iii) $r=3(c-1)/2$.
Thus the possibilities for $(c,d,k,r,\ell)$ are $(9,5,12,12,3)$, $(15,8,18,21,3)$, $(33,17,36,48,3)$.
\end{enumerate}

\medskip\noindent
\fbox{The case $\la=4$}\quad
Lemma \ref{le}(viii) yields  five possibilities for $(\ell,x)$, namely  $(2,1)$,  $(2,2)$, $(2,3)$, $(3,1)$, $(4,1)$.
We now split the study into these  cases.
\begin{enumerate}[(i)]
\item $(\ell,x)=(2,1)$. By Lemma \ref{le}(vii) $k\mid 48$ and by Lemma \ref{le}(x) $k\geq 4$. Moreover $k$ is even by Lemma \ref{le}(ii).
Thus $k\in\{4, 6, 8, 12, 16, 24, 48\}$. By Lemma \ref{le}(iv) $c=k-2$, and Lemma \ref{le}(vi) yields $k\mid 8c(c-1)$ which is satisfied in each case. Combining Lemma \ref{le}(iii) and (v) yields $d=c$. By Lemma \ref{le}(iii) $r=4(c-1)$.
Thus the possibilities for $(c,d,k,r,\ell)$ are  $(2,2,4,4,2)$, $(4,4,6,12,2)$, $(6,6,8,20,2)$, $(10,10,12,36,2)$, $(14,14,16,42,2)$, $(22,22,24,84,2)$, $(46,46,48,180,2)$. However, in the first case, $k=v$ so this design is trivial.

\item $(\ell,x)=(2,2)$. 
By Lemma \ref{le}(vii) $k\mid 96$ and by Lemma \ref{le}(x) $k\geq 8$. 
Thus $k\in\{8, 12, 16, 24, 32, 48, 96\}$. By Lemma \ref{le}(iv) $c=k/2-1$, and Lemma \ref{le}(vi) yields $k\mid 12 c(c-1)$ which is not satisfied for $k=16$, $k=48$ and $k=96$. Combining Lemma \ref{le}(iii) and (v) yields $d=2c-1$. By Lemma \ref{le}(iii) $r=4(c-1)$.
Thus the possibilities for $(c,d,k,r,\ell)$ are $(3,5,8,8,2)$, $(5,9,12,16,2)$, $(11,21,24,40,2)$, $(15,29,32,56,2)$.

\item $(\ell,x)=(2,3)$. 
By Lemma \ref{le}(vii) $k\mid 160$ and by Lemma \ref{le}(x) $k\geq 20$. 
Thus $k\in\{20, 32, 40, 80, 160\}$. By Lemma \ref{le}(iv) $c=(k-2)/3$ so $k$ cannot be $40$ nor $160$.  Lemma \ref{le}(vi) yields $k\mid 16 c(c-1)$ which is satisfied in each remaining case. Combining Lemma \ref{le}(iii) and (v) yields $d=3c-2$. By Lemma \ref{le}(iii) $r=4(c-1)$.
Thus the possibilities for $(c,d,k,r,\ell)$ are $(6,16,20,20,2)$, $(10,28,32,36,2)$, $(26,76,80,100,2)$.

\item $(\ell,x)=(3,1)$. 
By Lemma \ref{le}(vii) $k\mid 96$ and by Lemma \ref{le}(x) $k\geq 8$. Moreover $k$ is divisible by $3$ by Lemma \ref{le}(ii).
Thus $k\in\{ 12,  24,  48, 96\}$. By Lemma \ref{le}(iv) $c=k-3$, and Lemma \ref{le}(vi) yields $k\mid c(c-1)(k-2)$ and thus also $k\mid2 c(c-1)$, which is  not satisfied for  $k=48$ and $k=96$. Combining Lemma \ref{le}(iii) and (v) yields $d=(c+1)/2$.  By Lemma \ref{le}(iii) $r=2(c-1)$.
Thus the possibilities for $(c,d,k,r,\ell)$ are $(9,5,12,16,3)$, $(21,11,24,40,3)$.

\item $(\ell,x)=(4,1)$. 
By Lemma \ref{le}(vii) $k\mid 160$ and by Lemma \ref{le}(x) $k\geq 20$. 
Thus $k\in\{20, 32, 40, 80, 160\}$. By Lemma \ref{le}(iv) $c=k-4$. Combining Lemma \ref{le}(iii) and (v) yields $d=(c+2)/3=(k-2)/3$ so $k$ cannot be $40$ nor $160$. 
 Lemma \ref{le}(vi) yields $k\mid \frac{4 c(c-1)(c+2)}{9}$ which is  satisfied in each remaining case. By Lemma \ref{le}(iii) $r=4(c-1)/3$.
Thus the possibilities for $(c,d,k,r,\ell)$ are $(16,6,20,20,4)$, $(28,10,32,36,4)$, $(76,26,80,100,4)$.
\end{enumerate}
\qed
 
 However not all numerically feasible tuples listed above lead to an example, as we will see in Section \ref{sec:class}.

\section{A flag-regular, point-imprimitive design on \texorpdfstring{$36$}{36} points}\label{sec:36}

In this section we construct a flag-transitive, point-imprimitive design corresponding to the numerically feasible parameter tuple: 

\begin{equation}\label{eq:36tuple}
(\la,v,k,r,c,d,\ell)=(4,36,8,20,6,6,2) 
\end{equation}
from Proposition~\ref{la=3}. We also prove in Proposition~\ref{lem:36unique} that, up to isomorphism, this example is the unique flag-transitive design with the parameter set \eqref{eq:36tuple}.  Moreover, the design  satisfies equality in the bound $v\leq (k-2)^2$ of \cite[Proposition 4.1]{CP93}.
This design is hard to find in the literature, and we have sought advice from colleagues Alfred Wasserman, Patric \"{O}sterg\aa rd and Charles Colbourn. Collectively we were unable to find it. The Handbook of Combinatorial Designs \cite{Handbook} mentions references for two designs with parameters $(\la,v,k,r)=(4,36,8,20)$. Firstly in \cite{Abel} an example with these parameters is given with `repeated blocks', and secondly a construction in \cite{VT} produces an example which we were able to construct computationally; we found that its automorphism group has order $2$. Thus the design we present is neither of the ones listed in \cite{Handbook}. After completing the first draft of this paper we became aware of a new paper of Zhang and Zhou in which this new design occurs \cite[Theorem 1.3]{ZZ}. We make some comments about their work in Remark~\ref{rem:shenglin} below.  

We give several constructions for this design based on the symmetric group $\Sym_6$ of degree $6$. The first description gives sufficient information for the design to be constructed computationally, see also Remark~\ref{rem:con1}. It is based on the transitive permutation representation of $\Sym_6$ on $36$ points, and relies on an explicit description of an outer automorphism $\sigma$ of $\Sym_6$, namely $\sigma$ is determined by its action on a standard generating set for $\Sym_6$ as follows: 
\begin{align}\label{eq:sigma}
    (1,2)^\sigma&=(1,4)(2,6)(3,5)\quad\text{and}\quad (1,2,3,4,5,6)^\sigma=(1,3)(2,6,5).
\end{align}

\begin{construction}\label{con1}
Let $\P=\{(i,j) \mid\ 1\leq i, j\leq 6\}$, and let $G=\{(g,g^\sigma)|g\in\Sym_6\}$ acting coordinate-wise on $\P$. Let
$$
B=\{(1,1),(1,2),(2,1),(2,3),(3,2),(3,4),(4,3),(4,4)\}, 
$$
and let $\mathcal{B}=\{B^g \mid\, g\in G\}$, the $G$-orbit of $B$ under $G$, and define the design $\mathcal{D}= (\P, \mathcal{B})$.
\end{construction}

\begin{lemma}\label{lem:con1}
The design $\mathcal{D}= (\P, \mathcal{B})$ of Construction~$\ref{con1}$ is a $2-(36,8,4)$ design
with full automorphism group $G\cong\Sym_6$ acting flag-regularly and point imprimitively. Moreover, $G$ leaves invariant two nontrivial point-partitions, each with $d=6$ parts of size $c=6$, namely the `rows' and the `columns' of the square array $\P$. 
\end{lemma}

\proof
By definition, $G$ is admitted as an automorphism group of $\mathcal{D}$, and leaves invariant the two nontrivial point-partitions formed by the rows and the columns of $\P$. Also $\mathcal{D}$ has $v=36$ points and block size $k=8$. A computation using Magma~\cite{magma} yields that $\mathcal{D}$ is a $2$-design with $\lambda=4$ and that $G$ is the full automorphism group. 
\qed

\begin{remark}\label{rem:con1}
Computationally, using Magma~\cite{magma}, the design $\mathcal{D}$ of Construction~\ref{con1} can be constructed using the unique smallest block-transitive subgroup of automorphisms (which we note  is not flag-transitive on $\D$), namely the index $2$ subgroup $H=\Alt_6$ of $\Aut(\D)=\Sym_6$. The group $H$ can be constructed up to conjugacy in $\Sym_{36}$, using Magma, as $H =$  {\tt TransitiveGroup(36, 555)}. Then the block-set of $\D$ can be  constructed as the set of images of the $8$-element subset $B=\{ 1, 2, 7, 8, 22, 23, 25, 26\}$ of $\{ 1, 2,\dots, 36\}$ under the action of $H$. See also Table~\ref{table1}.

\end{remark}


\begin{remark}\label{rem:con2}
Construction~\ref{con1} gives some insight into the set of points and the group action on them, but  provides little understanding of the nature of the blocks. We now give a different construction for (a design isomorphic to) $\mathcal{D}$ which gives a better understanding of the blocks in terms of the standard action of $\Sym_6$ of degree $6$. 

 For this description we note that $G=\Sym_6$ has a unique conjugacy class of subgroups of index $36$, namely the class of Frobenius subgroups $F\cong F_{20}$ of order $20$. This means that we may identify the point set $\P$ with $\{ F^g\mid\,g\in G\}$,  where $G$ acts by conjugation. 

 Now $G$ has two conjugacy classes of subgroups of index $6$, corresponding to $\Sym_5$ and $\PGL_2(5)$, which are interchanged by the outer automorphism $\sigma$ given in \eqref{eq:sigma}.  Each Frobenius group in $\P$ is contained in a unique subgroup from each of these classes, giving two distinct $G$-invariant partitions of $\P$, each with $d=6$ parts of size $6$, see \cite[Lemma 2.14]{PS}.  

 For the construction, we use one of these partitions:   let $X=\{1,2,3,4,5,6\}$ be the set on which $G$ acts naturally, and note that each $F\in\P$ fixes a unique element of $X$. For each $x\in X$, let $\Delta_x$ be the set of six Frobenius groups in $\P$ which fix $x$, so that $\mathcal{C}=\{\Delta_x \mid\,x\in X\}$ is one of the $G$-invariant partitions described in the previous paragraph. (The second point-partition is based in a similar fashion on the set $Y$ of six transitive subgroups $\PGL_2(5)$.)  

 The blocks of the design are labelled by triples of the form $(x, x’, \pi)$, where $x, x’$ are distinct elements of $X$ and $\pi$ is a `bisection' of $X\setminus\{x,x'\}$, that is a partition with two parts of size $2$. For each pair $(x,x')$ there are three choices for $\pi$ and hence there are $6\times 5\times 3 = 90$ triples, hence $90$ blocks. 

 We need to define the $8$-subset of $\P$ forming 
the block $B=B(x, x', \pi)$. We note that for each of the four elements $z\in X\setminus\{x,x’\}$,  there are six Frobenius groups in $\Delta_z$, giving a set $\mathcal{P}(\pi)$ of $4\times 6=24$ points of $\P$. 
Using Magma, we find that the setwise stabiliser $H:=G_{\mathcal{P}(\pi)}$ has three orbits of length $8$ in $\mathcal{P}(\pi)$. One of these orbits is the block $B(x,x',\pi)$, and one of the other orbits is the block $B'=B(x',x,\pi)$ (which has the same stabiliser $G_{B'}=G_B=G_{\mathcal{P}(\pi)}$). The normaliser $N_G(H)$ interchanges the triples $(x, x’, \pi)$ and $(x', x, \pi)$, and we find (with Magma) that $N_G(H)$ interchanges two of the $H$-orbits in $\mathcal{P}(\pi)$ and leaves the third invariant. We choose one of the two $H$-orbits  
moved by $N_G(H)$  and call it $B(x,x',\pi)$, and we take  the block set $\mathcal{B}$ of $\mathcal{D}$ to be the set of $G$-images of this $8$-subset $B(x,x',\pi)$. Thus $\mathcal{D}=(\P,\mathcal{B})$ is well defined.


 It is clear from the construction that $\mathcal{D}$ is a flag-transitive point-imprimitive $1$-design. It may be checked using Magma that $\mathcal{D}$ is in fact a $2-(36,8,4)$ design. 
We note, finally, that the outer automorphism $\sigma$ in \eqref{eq:sigma} is not an automorphism of  $\mathcal{D}$, but rather $\sigma$ maps $\mathcal{B}$ to a different collection of ninety $8$-element subsets of $\P$ which forms a design isomorphic to $\mathcal{D}$. 
\end{remark}

Finally we prove that there is, up to isomorphism, a  unique flag-transitive point-imprim\-it\-ive design with parameters as in \eqref{eq:36tuple}, and it follows from this that the designs in Construction~\ref{con1} and Remark~\ref{rem:con2} are isomorphic.

\begin{proposition}\label{lem:36unique} 
Up to isomorphism, the design $\mathcal{D}$ in Construction~$\ref{con1}$  is the unique flag-transitive point-imprimitive design $\D=(\P,\B)$ with parameter set as in $\eqref{eq:36tuple}$.
\end{proposition}

\proof
Suppose that $\cal D=(\P,\B)$ has parameters  $(\la,v,k,r,c,d,\ell)=(4,36,8,20,6,6,2)$ and admits a flag-transitive point-imprimitive group $G$. By Proposition~\ref{la=3},  this is the only 
parameter tuple with $(\la,v,k)=(4,36,8)$, and hence each nontrivial $G$-invariant point-partition has $6$ parts of size $6$.  Let $\mathcal{C}=\{\Delta_1,\dots,\Delta_6\}$ be a $G$-invariant partition of $\P$ with each $|\Delta_i|=6$, and let $D=G^\mathcal{C}$ and $L=(G_{\Delta_1})^{\Delta_1}$. We may assume that $G\leq L\wr D\leq \Sym_6\wr \Sym_6$, by \cite[Theorem 5.5]{PS}.  Moreover both $L, D$ are primitive of degree $6$, as otherwise there would be another parameter set in Proposition \ref{la=3} for $(\la,v,k)=(4,36,8)$. This implies that $D$ and $L$ are $2$-transitive, and each has socle $\PSL_2(5)$ or $\Alt_6$. Thus, for distinct $i, j$, $G_{\{\Delta_i, \Delta_j\}}$ has index $\binom{6}{2}=15$ in $G$. 
 
 Now each block $B\in\B$  meets each of four parts $\Delta_i\in\mathcal{C}$ in $\ell=2$ points and is disjoint from the remaining two parts, and by  Lemma \ref{condition 1},  $b=|\B|=vr/k=90$. 
Thus $G_B$ has index $90$ in $G$ and fixes setwise the two parts, say $\{\Delta_i, \Delta_j\}$, which $B$ intersects trivially. Hence $G_B<G_{\{\Delta_i, \Delta_j\}}<G$ and $|G_{\{\Delta_i,\Delta_j\}}:G_B|=90/15=6$.  In particular $G_B$ contains all Sylow $5$-subgroups of $G_{\{\Delta_i, \Delta_j\}}$. Let $P$ be a Sylow $5$-subgroup of $G_{\{\Delta_i, \Delta_j\}}$, so $P\leq G_B$. Since the group induced by $G_{\{\Delta_i, \Delta_j\}}$ on $\mathcal{C}$ is a subgroup of $\Sym_2\times\Sym_4$, the order of which is not divisible by $5$, it follows that $P$ is contained in $K=G_{(\mathcal{C})}$, the kernel of the $G$-action on $\mathcal{C}$. Note that $K\leq G_{\{\Delta_i, \Delta_j\}}$, so a  Sylow $5$-subgroup of $G_{\{\Delta_i, \Delta_j\}}$ is also a  Sylow $5$-subgroup of $K$ and vice-versa.

Suppose that $K\ne 1$. Since $K$ is normal in $G$, its orbits on points all have the same size. In particular, for each $\Delta\in\mathcal{C}$, $K^{\Delta}$ is a nontrivial normal subgroup of the primitive group $L=G_{\Delta}^{\Delta}$, and hence $K^{\Delta}$ contains the socle of $L$. Since this socle is $\PSL_2(5)$ or $\Alt_6$, it follows that $5$ divides $|K^\Delta|$, and hence for some choice of Sylow $5$-subgroup $P$ of $K$ we have $P^\Delta\ne 1$. Since all Sylow $5$-subgroups of $K$ are conjugate in $K$, and since $K$ fixes each $\Delta\in\mathcal{C}$ setwise, this implies that,  for all $\Delta\in\mathcal{C}$, $P^\Delta\ne 1$ and has orbits of lengths $1, 5$ in $\Delta$. However, if $\Delta\not\in\{\Delta_i, \Delta_j\}$, then $P$ fixes setwise the $2$-subset $B\cap \Delta$ since $P\leq G_B$, which is a contradiction.  

Hence $K=1$, so $G\cong G^\mathcal{C}\leq \Sym_6$.  However also $|G|$ is divisible by the number of flags, which is $90\times 8 = |\Sym_6|$. Hence $G\cong \Sym_6$  and $G$ is regular on flags. Now $\Sym_6$ has a unique conjugacy class of subgroups of index $36$, namely the class of Frobenius groups $F_{20}$, and each such subgroup is contained in two distinct subgroups of index $6$ in $G$. Hence $G$ leaves invariant two distinct point-partitions with six parts of size six. This  unique transitive action of $\Sym_6$ of degree $36$ can be found as   {\tt TransitiveGroup(36,1252)} in Magma.
We checked with Magma~\cite{magma} that, up to isomorphism, the design $\D$ is as in Construction~\ref{con1}. To do this we  searched in  {\tt TransitiveGroup(36,1252)} for orbits of size $b=90$ on the set of $20250$ $8$-subsets which have $2$ points in each of four parts of both nontrivial invariant partitions. There are five such orbits, only two of which yield  $2$-designs, and the designs are isomorphic.
\qed

\begin{remark}\label{rem:shenglin}
After completing our work, the paper \cite{ZZ} of Zhang and Zhou appeared, in which this unique flag-transitive $2-(36,8,4)$ design arises in their classification \cite[Theorem 1.3]{ZZ} of $2$-designs with $\lambda\leq 4$ admitting an automorphism group which is flag-transitive, and acts imprimitively and quasiprimitively on points. We make a few comments about this result. The analysis of this case in their proof  \cite[pp. 431-432]{ZZ} does not give much detail. It is likely that their proof relies on extensive computation, but no references are made to this. Their proof proceeds by asserting (i) that there are seven conjugacy classes of subgroups of $G=S_6$ of index $90$ (the number of blocks), (ii) that five of these classes consist of subgroups $H$ which have an orbit $B$ on points of size $8$ such that the number of $G$-images of $B$ is $90$, and (iii) that ``it is easy to see that there are only two conjugacy classes of subgroups'' such that the set of $G$-images of the $H$-orbit $B$ is the block set of a $2$-design. 
They also state that ``it is not hard to check that'' the two designs obtained are isomorphic. No details are given, either of the construction or of any computations they may have carried out to justify the assertions.  

In  \cite[Theorem 1.3(ii)]{ZZ} it is stated that there is a unique non-symmetric $2-(36,8,4)$ design. Our Proposition~\ref{lem:36unique} proves that there is a unique such design which is flag-transitive and point-imprimitive. However, as discussed at the beginning of this section, there exists at least one other $2$-design with these parameters obtained from a construction in \cite{VT}. 
\end{remark}

\section{Classification for $v$ less than $100$ and $\lambda$ at most $4$ } \label{sec:class}
   
In Theorem \ref{class} we obtain a full classification of the designs for $v<100$ and $\la\leq 4$. This result will follow from the following citations and new results, and from the previous section. Up to isomorphism, we determine exactly eleven designs (two of the four designs with $v=96$ admit two distinct $(c,d)$ possibilities), and these correspond to exactly seven of the twelve numerically feasible tuples $(\lambda, v, k, r, c, d, \ell)$ in Proposition~\ref{la=3} with $\lambda\leq 4, v<100$. 

\begin{table}[ht]
\caption{\label{table1}Construction of all examples with $\la\leq 4$ and $v<100$, note the groups $H$ listed are block-transitive and as small as possible}
\begin{tabular}{|cccc|p{6.5cm}|p{2.3cm}|p{1.8cm}|cc|}
\hline
$\la$&$v$&$k$&$r$& Block $B$ & Group $H$ & $|\Aut(\cal D)|$&$c$&$d$\\
\hline
$2$&$16$&$6$&$6$& $\{ 1, 2, 7, 9, 12, 13\}$&{\tt TG(16,3)}&$11520$&$4$&$4$\\  
\cline{5-9}
&&&&$\{ 1, 2, 3, 7, 11, 15\}$ & {\tt TG(16,5)}&$768$&$4$&$4$\\
\hline
$3$&$45$&$12$&$12$&$\{1, 2, 3, 4, 8, 10, 20, 22, 25, 34, 39, 41\}$&{\tt TG(45,63)} &$19440$&$9$&$5$\\ 
\hline
$4$&$15$&$8$&$8$&$\{ 1, 2, 3, 4, 8, 11, 12, 14 \}$&{\tt TG(15,1)}&$20160$&$3$&$5$\\ 
\hline
$4$&$16$&$6$&$12$&$\{1,2,3,6,11,13\}$ & {\tt TG(16,27)}&$6144$&$4$&$4$\\ 
\cline{5-9}
&&&&$\{1,2, 5, 7, 13, 16\}$ & {\tt TG(16,46)}&$1920$&$4$&$4$\\ 
\hline
$4$&$36$&$8$&$20$&$\{ 1, 2, 7, 8, 22, 23, 25, 26\}$&{\tt TG(36,555)} &$720$&$6$&$6$\\
\hline
$4$&$96$&$20$&$20$&$\{ 1, 2, 3, 11, 12, 141, 17, 29, 31, 32, $&$H_1$ &$552960$&$16$&$6$\\
&&&&$36, 41, 47, 51, 57, 63, 68, 85, 87, 93\}$&&&&\\
\cline{5-9}
&&&&$\{ 1, 2, 3, 11, 14, 17, 24, 29, 31, 35, 43,$&$H_1$ &$184320$&$16$&$6$\\
&&&&$ 44, 48, 56, 64, 65, 69, 90, 95, 96\}$&&&$6$&$16$\\
\cline{5-9}
&&&&$\{  1, 2, 3, 5, 11, 15, 21, 22, 23, 27, 41, $&$H_2$ &$138240$&$16$&$6$\\
&&&&$ 43, 62, 68, 77, 80, 86, 90, 92, 95 \}$&&&&\\
\cline{5-9}
&&&&$\{1, 2, 3, 5, 11, 15, 21, 23, 31, 34, 46, $&$H_2$ &$7680$&$16$&$6$\\
&&&&$ 58, 66, 67, 69, 70, 72, 73, 89, 96\}$&&&$6$&$16$\\
\hline
\end{tabular}
\end{table}

\begin{remark}
In Table~\ref{table1} we list a block $B$ and group  $H$ which allow a quick construction of the corresponding design: namely $|H|$ is minimal  such that $H$ is block-transitive (not necessarily  flag-transitive) and hence the block-set is $B^H = \{ B^h\mid h\in H\}$. We also list the order of the full automorphism group. In most cases the group $H$ is described as {\tt TG(v,i)}, which is an abbreviation of the  name {\tt TransitiveGroup(v,i)} for the $i^{th}$ group of degree $v$ in the database of transitive groups of small degree in Magma~\cite{magma}. 
Since the transitive groups in Magma are only given for degrees up to $47$, we cannot describe $H$ in this way for the four designs with $v=96$ points.  In these cases we construct the designs using the method described in \cite{LPR09}: for each design we give generators for a group $H$  which is block-regular  (hence not flag-transitive). In all four cases the group $H$ is one of the following  two groups $H_1$, $H_2$.
\begin{align*}     
H_1=\langle & (1, 37, 2, 31)(3, 23, 8, 19)(4, 21, 9, 18)(5, 27, 7, 14)(6, 11, 10, 26)(12, 39, 24, 50)(13, 17, 25, 29)\\&(15, 40, 22, 42)(16, 38, 28, 20)(30, 33, 49, 45)(32, 55, 53, 65)(34, 47, 46, 35)(36, 43, 48, 41)\\&(44, 60, 62, 51)(52, 79, 70, 90)(54, 66, 63, 57)(56, 64, 58, 59)(61,
        87, 81, 94)(67, 91, 78, 95)\\&(68, 76, 80, 74)(69, 83, 89, 73)(71, 96, 75, 86)(72, 88, 82, 93)(77, 92, 84, 85),\\&
    (1, 62, 42)(2, 44, 40)(3, 29, 32)(4, 63, 28)(5, 59, 26)(6, 56, 27)(7, 64, 11)(8, 17, 53)(9, 54, 16)\\&(10, 58, 14)(12, 69, 81)(13, 23, 55)(15, 51, 37)(18, 38, 66)(19, 65, 25)(20, 57, 21)(22, 60, 31)\\&(24, 89, 61)(30, 95, 75)(33, 86, 67)(34, 68, 52)(35, 90, 74)(36, 92,
        72)(39, 94, 83)(41, 93, 77)\\&(43, 88, 84)(45, 96, 78)(46, 80, 70)(47, 79, 76)(48, 85, 82)(49, 91, 71)(50, 87, 73),\\&
    (1, 73, 2, 83)(3, 85, 8, 92)(4, 52, 9, 70)(5, 67, 7, 78)(6, 61, 10, 81)(11, 94, 26, 87)(12, 41, 24, 43)\\&(13, 80, 25, 68)(14, 91, 27, 95)(15, 71, 22, 75)(16, 93, 28, 88)(17, 76, 29, 74)(18, 79, 21, 90)\\&(19, 77, 23, 84)(20, 72, 38, 82)(30, 46, 49, 34)(31, 69, 37, 89)(32,
        64, 53, 59)(33, 47, 45, 35)\\&(36, 39, 48, 50)(40, 86, 42, 96)(44, 57, 62, 66)(51, 54, 60, 63)(55, 56, 65, 58) \rangle
\end{align*}
\begin{align*}
H_2=\langle&(1, 55, 18, 26)(2, 90, 24, 66)(3, 50, 38, 51)(4, 49, 5, 60)(6, 61, 35, 63)(7, 80, 10, 89)(8, 48, 47, 76)\\&(9, 83, 39, 84)(11, 96, 40, 86)(12, 53, 16, 57)(13, 93, 23, 94)(14, 65, 15, 44)(17, 85, 21, 87)\\&(19, 56, 20, 64)(22, 77, 34, 62)(25, 67, 33, 95)(27, 70, 58, 69)(28,
        92, 29, 68)(30, 73, 59, 72)\\&(31, 71, 43, 82)(32, 91, 42, 88)(36, 52, 37, 74)(41, 78, 46, 79)(45, 81, 54, 75),\\&
    (1, 75, 21)(2, 6, 79)(3, 48, 24)(4, 73, 40)(5, 69, 9)(7, 12, 81)(8, 61, 66)(10, 18, 71)(11, 37, 70)\\&(13, 19, 62)(14, 92, 23)(15, 77, 33)(16, 82, 17)(20, 68, 25)(22, 64, 95)(26, 80, 45)(27, 49, 96)\\&(28, 56, 94)(29, 44, 67)(30, 74, 86)(31, 53, 89)(32, 38, 78)(34, 65,
        93)(35, 76, 42)(36, 72, 39)\\&(41, 63, 91)(43, 55, 85)(46, 50, 90)(47, 51, 88)(52, 83, 58)(54, 57, 87)(59, 60, 84),\\&
    (1, 38)(2, 39)(3, 18)(4, 15)(5, 14)(6, 12)(7, 13)(8, 58)(9, 24)(10, 23)(11, 42)(16, 35)(17, 25)\\&(19, 36)(20, 37)(21, 33)(22, 54)(26, 60)(27, 47)(28, 31)(29, 43)(30, 41)(32, 40)(34, 45)(44, 51)\\&(46, 59)(48, 92)(49, 55)(50, 65)(52, 53)(56, 61)(57, 74)(62, 79)(63, 64)(66,
        89)(67, 96)(68, 76)\\&(69, 71)(70, 82)(72, 81)(73, 75)(77, 78)(80, 90)(83, 94)(84, 93)(85, 88)(86, 95)(87, 91)\rangle
        \end{align*}
\end{remark}

\begin{remark}\label{rem:t2}
In Table \ref{table2}, we give additional information about the groups of these designs, obtained using Magma~\cite{magma}. In column $\Aut(\cal D)$, we list the full automorphism group of the design: for $v\ne 15, 96$, we give the group in the form  {\tt TransitiveGroup(v,i)} as well as its structure; for $v=15$, $\Aut(\cal D)$ is a well known group given in its standard action; while for $v=96$ (which is not covered by the database in \cite{magma}) we give the structure of the group, but note that the full automorphism group can be easily found in Magma~\cite{magma} by constructing the design using data in Table \ref{table1} and calling for its full automorphism group. In column \emph{Largest FT imp.}, we list, up to conjugacy, the largest flag-transitive subgroup of $\Aut(\cal D)$ which preserves a partition with $d$ parts of size $c$. If $\Aut(\cal D)$ itself preserves such a partition, then we just write $\Aut(\cal D)$. We draw attention to exceptional behaviour for two of the designs with  $v=96$, namely the second and fourth designs in the last block of Table~\ref{table2}, which are the designs numbered $2$ and $4$, respectively, in \cite[Table 1]{LPR09}. For these designs, $\Aut(\cal D)$ preserves a partition with $6$ blocks of size $16$ but not a partition with $16$ blocks of size $6$, while a proper flag-transitive subgroup preserves both.  In column \emph{Smallest FT imp.}, we list, up to conjugacy, the flag-transitive subgroups of $\Aut(\cal D)$ which preserve a partition with $d$ parts of size $c$ and are of smallest size. Note there is not always a unique such subgroup, as shown in the table.
\end{remark}

\begin{table}[ht]
\caption{\label{table2}Full automorphism groups, largest flag-transitive subgroup preserving the partition (with the given $c$ and $d$), smallest flag-transitive subgroups preserving this partition }
\begin{tabular}{|p{.3cm}p{.3cm}p{.3cm}p{.3cm}p{.3cm}p{.3cm}|p{3.8cm}|p{3cm}|p{4.9cm}|}
\hline
$\la$&$v$&$k$&$r$&$c$&$d$& $\Aut(\cal D)$ & Largest FT imp.& Smallest FT imp.\\
\hline
$2$&$16$&$6$&$6$&$4$&$4$&{\tt TG(16,1753)}$\cong
2^4.\Sym_6$&{\tt TG(16,1063)}$\cong  $&{\tt TG(16,183)}$\cong  2^4. 6$,  \\
 &&&&&&&$2^4.(\Sym_2\wr\Sym_3)\cong  $&{\tt TG(16,184)}$\cong
 2^3.\Alt_4$, \\
 &&&&&&&$2^5.\Sym_4$&{\tt TG(16,185)}$\cong
 2^3.\Alt_4$,\\
  &&&&&&&& {\tt TG(16,194)}$\cong
 2^4.\Sym_3$, \\
 &&&&&&&&{\tt TG(16,195)}$\cong
 2^2.\Sym_4$\\
\cline{7-9}
&&&&&&{\tt TG(16,1073)}$\cong  2^5.\Sym_4$&$\Aut(\cal D)$&$\Aut(\cal D)$\\ \hline
$3$&$45$&$12$&$12$&$9$&$5$&{\tt TG(45,628)}$\cong  3^4. 2.\Sym_5$&$\Aut(\cal D)$&{\tt TG(45,314)}$\cong  3^4. 2.\AGL_1(5)$\\ \hline
$4$&$15$&$8$&$8$&$3$&$5$&$\PSL_4(2)$&{\tt TG(15,21)}$\cong  3.\PGGL_2(4) \cong (\Alt_5\times 3).2$&$\PGGL_2(4)\cong\Sym_5$\\ \hline
$4$&$16$&$6$&$12$&$4$&$4$&{\tt TG(16,1690)}$\cong  2^8.\Sym_4$&$\Aut(\cal D)$&  {\tt TG(16,419)}$\cong
 2^4.\Alt_4$,  \\ 
 &&&&&&&&{\tt TG(16,420)}$\cong  2^4. \Alt_4$ (two classes), \\
  &&&&&&&&{\tt TG(16,430)}$\cong
 2^3.\Sym_4$, \\
   &&&&&&&&{\tt TG(16,433)}$\cong
 2^4.((2\times 3).2)$\\
\cline{7-9}
&&&&&&{\tt TG(16,1329)}$\cong  2^4.\Sym_5$&{\tt TG(16,776)}$\cong  2^4.\Sym_4$&{\tt TG(16,776)}$\cong  2^4.\Sym_4$\\ \hline
$4$&$36$&$8$&$20$&$6$&$6$&{\tt TG(36,1252)}$\cong \Sym_6$&$\Aut(\cal D)$&$\Aut(\cal D)$\\ \hline
$4$&$96$&$20$&$20$&$6$&$16$&$  2^8.\Sym_6$&$ 2^4.\Sym_6$&$ 2^4.\Sym_5$\\
\cline{7-9}
&&&&&&$ 2^6.\Sym_5$&$ 2^4.\Sym_5$&$ 2^4.\Sym_5$\\
\hline
$4$&$96$&$20$&$20$&$16$&$6$&$  2^8.(( 3\times\Alt_6). 2)$&$\Aut(\cal D)$&$  2^8.\Alt_5$\\
\cline{7-9}
&&&&&&$ 2^8.\Sym_6$&$\Aut(\cal D)$&$  2^4.\Sym_5$\\
\cline{7-9}
&&&&&&$  2^6.(( 3.\Alt_6). 2)$&$\Aut(\cal D)$& $ 2^5.\Sym_5$\\
\cline{7-9}
&&&&&&$ 2^6.\Sym_5$&$\Aut(\cal D)$&$  2^4.\Sym_5$ and $  2^5.\Alt_5$\\
\hline
\end{tabular}
\end{table}

Suppose that $\D$ is   a non-trivial $2$-$(v,k,\la)$  design  admitting 
a flag-transitive point-imprimitive group $G$ of automorphisms, with $\lambda\leq4$ and $v<100$. 
As mentioned earlier $\la\geq 2$, and the classification for $\la=2$ is given in \cite{DLPX}. We use the work in \cite{DLPX} to describe the designs for $\la=2$. 
For $\la=3,4$,
 the tuple $(\la,v,k,r,c,d,\ell)$ must be   numerically feasible and so appears in the table in Proposition~\ref{la=3}.  We consider the possibilities for $\lambda$ separately. 

\medskip\noindent
\fbox{The case $\la=2$}

This case was analysed in  \cite[Theorem 1.1]{DLPX} where it was shown that that the only possible tuple is $(v,k,r,c,d,\ell)=(16,6,6,4,4,2)$, yielding two non-isomorphic examples, as in  lines 1--2 of Table~\ref{table1}. These two designs were first constructed by Hussain \cite{Huss} in 1945.
The first example has a flag-regular imprimitive group, the full automorphism group is {\tt TransitiveGroup(16,1753)} of order 11520 (which is primitive), and contains flag-transitive imprimitive subgroups, the largest, isomorphic to {\tt TransitiveGroup(16,1063)},  having order 768. 
The second  example has  full automorphism group {\tt TransitiveGroup(16,1073)}, which has order 768 and is imprimitive. This is also the unique flag-transitive subgroup of automorphisms.

\bigskip
\newpage\noindent
\fbox{The case $\la=3$}

In \cite{ZZ18}, Zhan and Zhou classify the flag-transitive imprimitive $2$-designs with $k=6$. It turns out they all have $v=16$, moreover there are two with $\la=2$ (as mentioned above), none with $\la=3$ and two with $\la=4$ (see below). 
So the tuple $(v,k,r,c,d,\ell)=(16,6,9,4,4,2)$ is not possible.

By \cite[Corollary 1.2]{P07}, there is a unique flag-transitive, point-imprimitive  $2-(45, 12, 3)$ design. It has (up to isomorphism) automorphism group {\tt TransitiveGroup(45,628)} which is imprimitive
preserving a partition with $d=5$ parts of size $c=9$, as in line 3 of Table~\ref{table1}.  The smallest imprimitive flag-transitive subgroup is isomorphic to {\tt TransitiveGroup(45,314)}.
Moreover by \cite[Proposition 5.1]{P07}, the tuple $(v,k,r,c,d,\ell)=(45,12,12,5,9,2)$ is not 
possible.  

\medskip\noindent
\fbox{The case $\la=4$}

The projective example described by Huw Davies in \cite{Dav1}, and mentioned in the introduction, (whose blocks are the hyperplane complements), provides an example in the case $n=3$ for which $(v,k,r,c,d,\ell)=(15,8,8,3,5,2)$, and in fact, by \cite[Proposition 1.5, see also Section 4.1]{PZ06}, it is the unique example, up to isomorphism. Its full automorphism group $\PSL(4,2)$ is point-primitive, while the largest flag-transitive point-imprimitive subgroup is isomorphic to {\tt TransitiveGroup(15,21)} of order $360$, and contains $\PGGL(2,4)\cong \Sym_5$ (which is regular on flags).

By \cite[Main Theorem and Table 3]{ZZ18}, the tuple  $(v,k,r,c,d,\ell)=(16,6,12,4,4,2)$ admits exactly two examples. 
For the first example,   {\tt TransitiveGroup(16,1690)} is the full automorphism group, of order $6144$; it is point-imprimitive and has a subgroup that is regular on flags.   
The second example  has full automorphism  group {\tt TransitiveGroup(16,1329)} of order $1920$, which is point-primitive; it has a unique flag-transitive imprimitive subgroup, and this subgroup is isomorphic to {\tt TransitiveGroup(16,776)} with flag stabiliser of order~$2$. 

In \cite[Theorem 1.1]{LPR09} it is shown that, up to isomorphism, there are exactly four flag-transitive $2-(96,20,4)$ designs, and for each of them the full automorphism group preserves a point-partition with $d=6$ parts of size $c=16$, see \cite[Subsections 1.2 and 3.1]{LPR09}. Thus the parameter tuple $(v,k,r,c,d,\ell)=(96,20,20,16,6,4)$ admits four examples. The full automorphism groups of these four designs have orders $552960$, $184320$, $138240$ and $7680$ respectively, and all flag-transitive subgroups of them have been determined, see \cite[Section 5]{LPR09}.  Each of the flag-transitive subgroups is listed in \cite[Table 2]{LPR09} and is of the form $C_2^a\rtimes H$ where $a\in\{4,5,6,8\}$. Assume $\cal D$ is one of these four designs where the  tuple $(v,k,r,c,d,\ell)=(96,20,20,6,16,2)$ is realised by a flag-transitive subgroup of automorphisms. Then, by \cite[Lemma 3.1]{LPR09}, $\Aut(\cal D)$ has a flag-transitive subgroup of the form $C_2^4\rtimes H$ (and a block-transitive subgroup $C_2^4\rtimes \Alt_5$), and hence by \cite[Table 2]{LPR09}, $\D$ is the design with $|\Aut(\cal D)|$ equal to either $7680$ or $184320$. 
We checked with Magma that in the first case the flag-transitive subgroup isomorphic to $C_2^4\rtimes \Sym_5$ admits two $G$-invariant partitions with $(c,d)=(6,16)$, while in the second case the flag-transitive subgroup isomorphic to $C_2^4\rtimes \Sym_5$ admits one $G$-invariant partition with $(c,d)=(6,16)$.

\medskip
By Proposition~\ref{lem:36unique}, there is, up to isomorphism, a unique flag-transitive point-imprimitive design with parameter set $(v,k,r,c,d,\ell)=(36,8,20,6,6,2)$, namely the design in Construction~\ref{con1}, and by Lemma~\ref{lem:con1} and Remark~\ref{rem:con1}, the entry in Table~\ref{table2} is valid. 

\medskip
The remaining parameter sets, both with  $v=45$, are dealt with in the following lemma.

\begin{proposition}\label{le:design45} 
There is no  flag-transitive, point-imprimitive $4$-design with parameter set  
$$
(v,k,r,c,d,\ell) \text{  equal to  } (45,12,16,9,5,3), \text{  or  } (45,12,16,5,9,2).
$$
\end{proposition}

\proof
Suppose that such a design $\D=(\P,\mathcal{B})$ exists, admitting a flag-transitive, 
point-imprimitive  automorphism group $G$.
Then $b = |\B| = vr/k=60$ and the number of flags is $f=bk=vr =720=2^4.3^2.5$. Thus $|G|=fz$ for some integer $z\geq1$.

Let $\mathcal{C}=\{\Delta_1,\dots,\Delta_d\}$ be a  $G$-invariant partition of the point-set $\P$ with each $|\Delta_i|=c$, 
where $(c,d)$ is either $(9,5)$ or $(5,9)$. Let $D=G^\mathcal{C}$, and $L=(G_{\Delta_1})^{\Delta_1}$, so by \cite[Theorem 5.5]{PS} we may assume that $G\leq L\wr D\leq \Sym_c\wr \Sym_d$, acting imprimitively. By Lemma \ref{la=3}, there  are no numerically feasible parameter sets with $(\la,v)=(4,45)$ and $c$ or $d$ equal to $3$, and hence both  $L$ and $D$ are primitive of degree $c$ and $d$ respectively. We note that each primitive group $X$ of degree $9$ has socle $T=C_3^2$ (affine type), or $PSL_2(8)$ or $\Alt_9$, and in the affine case, $T=O_3(X)$ is the largest normal $3$-subgroup of $X$, and is the group of translations (see, for example \cite[Theorem 3.15]{PS}). 

\medskip\noindent
\emph{Claim $1$:} If $(c,d)=(5,9)$, then there exists a second $G$-invariant partition of $\P$ 
with $5$ parts of size $9$, so that, without loss of generality we may assume that $(c,d)=(9,5)$.

\medskip\noindent
\emph{Proof of claim:} Suppose that $(c,d)=(5,9)$. In this case $\ell=2$, so a block meets each of six parts $\Delta_i$ in a $2$-subset and is disjoint from the remaining three parts. Now
\[
\#\{(B,\Delta) \mid B\in\B, \Delta\in\mathcal{C}, B\cap \Delta=\emptyset \} = b\times 3 = 9\times x,
\]
where $x$  is the number of blocks disjoint from a given class. Hence $x=b/3=20$, so each part $\Delta$ meets exactly $60-x=40$ blocks nontrivially. If $D$ has socle $PSL_2(8)$ or $\Alt_9$, then $|D|$, and hence also $|G|$, is divisible by $7$.  Since $b=60$, $7$ also divides $|(G_B)^\mathcal{C}|$, which is a contradiction since $(G_B)^\mathcal{C}$ fixes setwise the set of three parts disjoint from $B$. It follows that $C_3^2\unlhd D=G^\mathcal{C}\leq \AGL(2,3)$, and in particular $5$ does not divide $|G^\mathcal{C}|$. Since $b=60$ divides $|G|$, this implies that $5$ divides the order of  $K=G_{(\mathcal{C})}$, the kernel of the $G$-action on $\mathcal{C}$. In particular, $K\ne 1$. 
Since $K$ is normal in $G$, its orbits on points all have the same size, and hence the $K$-orbits are the parts $\Delta_i$ of $\mathcal{C}$, and since $K^{\Delta_i}$ is a transitive group of prime degree $5$, it is primitive. 

Next we show that $K$ acts faithfully on $\Delta_1$. If this is not the case then the kernel $K_{(\Delta_1)}$ of the action of $K$ on $\Delta_1$ is a nontrivial normal subgroup of $K$, and hence acts nontrivially on some part $\Delta\ne \Delta_1$. Thus $K_{(\Delta_1)}^\Delta$ is a nontrivial normal subgroup of the primitive group $K^\Delta$, and hence is transitive. This implies that a Sylow $5$-subgroup $P$ of $K_{(\Delta_1)}$ is nontrivial and acts transitively on $\Delta$. 
For each point-pair $\pi\subset \Delta_1$, there are exactly $\lambda=4$ blocks containing $\pi$, and since $P$ fixes $\pi$ (pointwise) it follows that $P$ fixes this set of 4 blocks setwise, and in fact $P$ fixes each of these four blocks (since $P$ is a $5$-group). Since this holds for all pairs $\pi\subset \Delta_1$, and since each block meeting $\Delta_1$ intersects it in two points, it follows that $P$ fixes setwise each of the 40 blocks which intersect $\Delta_1$ nontrivially. For any such block, say $B$, $B$ meets six parts in a two-subset, and each of these part-intersections with $B$ must be fixed setwise by $P$. It follows that $P$ must fix each of these six parts pointwise, and the same argument yields that $P$ fixes setwise each block which intersects any of these six parts nontrivially. Since each block intersects nontrivially with six of the nine parts of $\mathcal{C}$, it follows that $P$ fixes each block of $\cal B$ setwise. This contradicts the fact that $P$ is transitive on $\Delta$. Thus $K_{(\Delta_1)}=1$, and so
  $K\cong K^{\Delta_1}\leq \Sym_5$.   

Now we consider the map $\phi:G\to \Aut(K)$ induced by conjugation, and let $N = ker(\phi)=C_G(K)$. 
By the previous paragraph, $K\leq \Sym_5$, and in fact either $K=\Alt_5$ or $\Sym_5$, or the largest normal $5$-subgroup $O_5(K)$ of $K$  is isomorphic to $C_5$. In all cases $\Aut(K)$ is isomorphic to a subgroup of $\Sym_5$, and in particular $|\Aut(K)|$ is not divisible by $9$. Since $9$ divides $v$ and hence $|G|$, we conclude that $3$ divides $|N|$. Further, $N\cap K =C_G(K)\cap K=Z(K)$, and  either $Z(K)=1$ or $Z(K)=K=C_5<N$. In either case, $3$ divides the order of $N/(N\cap K)\cong N^\mathcal{C}$, which is a normal subgroup of the primitive group $G^\mathcal{C}=D$. Hence $N^\mathcal{C}$ contains the translation subgroup $T\cong C_3^2$ of $D\leq \AGL_2(3)$. 
Let $N_0$ be the (uniquely determined) subgroup of $N$ such that $N\cap K < N_0$ and $N_0^\mathcal{C}=T$,  and let $M$ be a Sylow $3$-subgroup of $N_0$. By definition $N_0$ is normal in $G$. Since $|N\cap K|=1$ or $5$, and since $N_0$ centralises $N\cap K$, it follows that $N_0\cong M\times (N\cap K)$ and in particular $M=O_3(N_0)\cong C_3^2$ is the unique Sylow $3$-subgroup of $N_0$. Thus $M$ is a characteristic subgroup of $N_0$, and hence is normal in $G$. Since $|M|=9$, the $M$-orbits in $\P$ form a $G$-invariant partition with $5$ parts of size $9$, which proves the claim. \qed

\medskip
Thus we may assume that $(c,d)=(9,5)$, so now  $\ell=3$, and  each block  meets each of four parts $\Delta_i$ in a $3$-subset and is disjoint from the remaining class. This time 
\[
\#\{(B,\Delta) \mid B\in\B, \Delta\in\mathcal{C}, B\cap \Delta=\emptyset \} = b\times 1 = 5\times x,
\]
where $x$  is the number of blocks disjoint from a given class, so $x=b/5=12$, and each part meets $60-x=48$ blocks nontrivially. 
Let $K=G_{(\mathcal{C})}$, the kernel of the $G$-action on $\mathcal{C}$. 

\medskip\noindent
\emph{Claim $2$:} $L=(G_{\Delta_1})^{\Delta_1}$ is of affine type, and $K^{\Delta_1}$ contains the translation group $O_3(L)\cong T$. Moreover,  for $Q=O_3(K)$, the largest normal $3$-subgroup of $K$, the $Q$-orbits and the $K$-orbits in $\P$ are the parts of $\mathcal{C}$, and $Q^{\Delta_i}\cong O_3(L)$ for each $i$. 

\medskip\noindent
\emph{Proof of claim:} 
Note that $G/K\cong D\leq \Sym_5$, and so $|G:K|$ is not divisible by $9$. Since $|G|=fz$ is divisible by $9$, it follows that $3$ divides $|K|$ and so  $K\ne 1$. Now $K$ is normal in $G$ and hence its orbits on points all have the same size. In particular $K^{\Delta_1}$ is nontrivial and normal in the primitive group $L$. Hence  $K^{\Delta_1}$ is transitive, and the $K$-orbits are the parts in $\mathcal{C}$. 
Let $B, B'$ be blocks which meet $\Delta_1$, say $\alpha\in B\cap \Delta_1, \alpha'\in B'\cap \Delta_1$. Since $G$ is flag-transitive, there exists $g\in G$ which maps the flag $(\alpha, B)$ to the flag $(\alpha', B')$ and hence $B^g=B'$ and $g$ fixes setwise the class $\Delta_1$ containing $\alpha$ and $\alpha'$. Thus $g\in G_{\Delta_1}$, and it follows that $G_{\Delta_1}$ is transitive on the set of $48$ blocks meeting $\Delta_1$ nontrivially. Thus $|G_{\Delta_1}:G_{\Delta_1, B}|=48$, and hence 
 $|L:G_{\Delta_1,B}^{\Delta_1}|$ divides $48$. If $L$ has socle $PSL_2(8)$ or $\Alt_9$ then $7$ divides $|L|$, and since   $|L:G_{\Delta_1,B}^{\Delta_1}|$ divides $48$, it follows that $7$ also divides $|G_{\Delta_1, B}^{\Delta_1}|$. This is a contradiction since $G_{\Delta_1,B}^{\Delta_1}$ leaves invariant the $3$-subset $B\cap \Delta_1$. Thus $L$ is of affine type, and hence $K^{\Delta_i}$ contains the translation group $O_3((G_{\Delta_i})^{\Delta_i})\cong T$, for each $i$. It follows that $Q=O_3(K)$ induces $T$ on each part $\Delta_i$ and hence the $Q$-orbits are the parts of $\mathcal{C}$. \qed


\medskip
We may therefore view each $\Delta_i$ as the affine plane $\AG_2(3)$.


\medskip\noindent
\emph{Claim $3$:} The $Q$-orbits in $\B$ have size $3$, and if $B\cap \Delta_i\ne\emptyset$, then  the $Q_B$-orbits in $\Delta_i$ form a parallel class of lines of the affine plane $\Delta_i$. Moreover, for each $i$, each line of the affine plane $\Delta_i$ occurs as the intersection with $\Delta_i$ of exactly four blocks, and each parallel class of lines of $\Delta_i$ corresponds to $12$ of these block-part intersections.  

\medskip\noindent
\emph{Proof of claim:} 
Since $Q$ is normal in $G$, the $Q$-orbits in $\B$ all have the same length, say $y$. So $y$ divides $b=60$, and $y$ is a power of $3$ since $y$ divides $|Q|$, whence $y=1$ or $3$. Since $Q^{\Delta_1}$ is transitive, it acts nontrivially on the blocks intersecting $\Delta_1$ in a $3$-subset, and hence $y=3$.

Since $G^\mathcal{C}$ is transitive, it is sufficient to prove the other assertions for $\Delta_1$. Let $B$ be a block such that  $B\cap \Delta_1\ne\emptyset$. Then $Q_B$ has index $3$ in $Q$, and as $Q^{\Delta_1}$ is the translation group by Claim 2, it follows that the $Q_B$-orbits in $\Delta_1$ form a parallel class of lines of the affine plane. 
 Let $\alpha\in \Delta_1$. Then $\alpha$ lies in $r=16$ blocks, and also $\alpha$ lies in four lines of the affine plane $\Delta_1$. For each of these lines $m$, and each point $\beta\in m\setminus\{\alpha\}$, the pair $\{\alpha,\beta\}$ lies in $\lambda=4$ blocks, each  intersecting $\Delta_1$ in the unique affine line $m$ containing $\{\alpha,\beta\}$.  Thus each of the affine lines on $\alpha$ occurs as the intersection with $\Delta_1$ of exactly four blocks. This is true for all points of $\Delta_1$, so each line of the affine plane $\Delta_1$ is the intersection with $\Delta_1$ of exactly four blocks. Moreover each parallel class of lines of $\Delta_1$ corresponds to $3\times 4$ block intersections with $\Delta_1$. \qed 

\medskip\noindent
\emph{Claim $4$:} $Q=T\cong C_3^2$ is faithful on each $\Delta_i\in\mathcal{C}$.  

\medskip\noindent
\emph{Proof of claim:} 
Since $Q^{\Delta_1}=T$ is the translation group, the subgroup $R=Q_{(\Delta_1)}$  fixes $\Delta_1$ pointwise, and is equal to $Q_\alpha$ for each $\alpha\in \Delta_1$. By Claim 3, for each of the 48 blocks $B$ that meet $\Delta_1$, the intersection $m=B\cap \Delta_1$ is a line of $\Delta_1$, and $Q_B=Q_m$ has index $3$ in $Q$. This implies that, for $\alpha\in m$, $Q_{m,\alpha}$ has index $3$ in $Q_m$ and hence index $9$ in $Q$, and we conclude that $Q_\alpha= Q_{m,\alpha} < Q_B < Q$.  Thus $R=Q_\alpha$ fixes each of the 48 blocks which meet $\Delta_1$.
Let $\Delta_i$ be one of the other three parts meeting such a block $B$. Then $R$ fixes $B\cap \Delta_i$ setwise, and hence each $R$-orbit in $\Delta_i$ is contained in a line of $\Delta_i$ parallel to $B\cap \Delta_i$.  By Claim 3, there are just 12 blocks which meet $\Delta_i$ in a line parallel to $B\cap \Delta_i$, while there are 48 blocks which meet  $\Delta_i$ nontrivially, and at most 12 of these are disjoint from $\Delta_1$. Hence there exists a block $B'$ which meets both $\Delta_1$ and $\Delta_i$ and is such that the line $B'\cap \Delta_i$ is not parallel to $B\cap \Delta_i$. We have shown that $R$ fixes each of the (non-parallel lines) $B\cap \Delta_i, B'\cap \Delta_i$ setwise, and hence $R$ fixes their intersection, which is a single point $\alpha'\in \Delta_i$. It follows that $R=Q_{\alpha'}$ and so $R$ fixes $\Delta_i$ pointwise, and hence fixes setwise every block meeting $\Delta_i$. Since this holds for each part $\Delta_i$ meeting $B$, it follows that $R$ fixes setwise every block that meets any of these four parts, and this implies that $R$ fixes every block of $\B$. Hence $R=1$, proving the claim. \qed

\medskip\noindent
\emph{Claim $5$:} $K\cong K^{\Delta_i}$  is faithful, for each $\Delta_i\in\mathcal{C}$.  

\medskip\noindent
\emph{Proof of claim:} As $G^\mathcal{C}$ is transitive it is sufficient to prove this for $\Delta_1$. Let $A=K_{(\Delta_1)}$, the pointwise stabiliser of $\Delta_1$ in $K$. By Claim 4, $A\cap Q=1$, and it follows that the normal subgroups $A, Q$ of $K$ centralise each other. Then for each $j$, $A^{\Delta_j}$ is contained in the centraliser of $Q^{\Delta_j}$ in $G_{\Delta_j}^{\Delta_j}\cong L$. 
Since $Q^{\Delta_j}$ is self-centralising in $G_{\Delta_j}^{\Delta_j}$ it follows that $A^{\Delta_j}\leq Q^{\Delta_j}$, and in particular $A^{\Delta_j}$ is a $3$-group.  Since $A$ is isomorphic to a subgroup  of $\prod_{j=1}^5 A^{\Delta_j}$, it follows that $A$ is a $3$-group. Thus $A\leq O_3(K)=Q$, and hence $A=1$. \qed 

\medskip
Since $D=G^\mathcal{C}\leq \Sym_5$,
it follows from Claims 2 and 5 that $|G|=|G^\mathcal{C}|.|K|$ divides $|\Sym_5|\times |\AGL_2(3)|=120\times 9\times 48$. Recall that $|G|=fz$ with $f=720=2^4.3^2.5$, the number of flags. Hence $z$ divides $72$.  
To complete this analysis we performed the following check computationally, using Magma:

\begin{itemize}
    \item We constructed the group $W=\AGL_2(3)\wr \Sym_5$ in its natural imprimitive permutation action on $\P$ of degree $45$ leaving invariant a partition $\mathcal{C}=\{\Delta_1,\dots,\Delta_5\}$ with each $|\Delta_i|=9$;
    \item for each subgroup $G$ of $W$ with order $fz$, for $z$ a divisor of $72$, we checked whether $G$ had an orbit $\B$ of size $b=60$ on $12$-subsets $B$ of $\P$ such that $B\cap \Delta_i$ is a line of the corresponding affine plane for exactly four parts $\Delta_i\in\mathcal{C}$; 
    \item for each such $G$ and $\B$, we checked whether $(\P, \B)$ was a $2$-design with $\la=4$.
\end{itemize}
This computer search yielded no $2$-designs.
 \qed

\section*{Acknowledgement}
The authors thank Charlie Colbourn, Patric \"{O}sterg\aa rd, and Alfred Wasserman for their advice about the $36$-point design.

\end{document}